\documentclass{llncs}
\usepackage{makeidx}  
\usepackage[pdftex]{graphicx}
\usepackage{amsmath}
\begin{document}
\frontmatter          
\pagestyle{headings}  
\addtocmark{Hamiltonian Mechanics} 
\mainmatter              
\title{Computing without a computer: a new approach for solving nonlinear differential equations}
\titlerunning{Computer Analogy}  
%
\author{Vladimir Aristov \and Andrey Stroganov }
\authorrunning{Vladimir Aristov, Andrey Stroganov} 
%
\tocauthor{Vladimir Aristov, Andrey Stroganov}
\institute{Dorodnicyn Computing Centre of Russian Academy of Sciences,\\
\email{aristov@ccas.ru}
}

\maketitle              

\begin{abstract}
The well-known Turing machine is an example of a theoretical digital computer, 
and it was the logical basis of constructing real electronic computers. In the 
present paper we propose an alternative, namely, by formalising arithmetic 
operations in the ordinary computing device, we attempt to go to the analytical 
procedure (for calculations). The method creates possibilities 
for solving nonlinear differential equations and systems. Our theoretical 
computer model requires retaining a finite number of terms to represent numbers, and
utilizes digit carry procedure. The solution is represented in the form of a segment of a 
series in the powers of the step size of the independent variable in the 
finite-difference scheme. The algorithm generates a schematic representation 
that approximates the convergent finite-difference scheme, which, in turn, 
approximates the equation under consideration. The use of probabilistic methods 
allows us to average the recurrent calculations and exclude intermediate levels 
of computation. All the stages of formalizing operations of the classical 
computer result in "the method of the computer analogy". The proposed method leads 
to an explicit analytical representation of the solution. We present the 
general features of the algorithm which are illustrated by an example of 
solutions for a system of nonlinear equations.  

\end{abstract}
\section{Introduction}
Analytical and numerical methods complement each other.  
Although most of nonlinear equations have no analytical solutions, 
yet numerical methods allow us to obtain solutions using the computer. 
It is important to emphasize 
that the computer is a material object (really, the electronic computer). 
In other words, an idealized mathematical object needs, in most cases, 
a materialized mathematical scheme. The question arose: can we formalize the 
operations of this material object itself  in such a way that the numerical 
(calculating) procedures will become analytical? 

Due to the efforts of some outstanding mathematicians such as Alan Turing, the 
computer theory has a solid base for constructing different types of computing 
machines (see [1]). The recent general tendency is as follows: most 
investigations aim at elaborating new computing devices, e.g. quantum computers 
[2] or developing analytical procedures (and languages) directed for improving 
computer software [3], i.e., their general purpose is to elaborate new logical 
computer schemes. Yet, we can discuss other proposals concerning analytical 
computing techniques (and we study not the logical but the arithmetical structure of 
the computer). 

We will consider a possibility to compute analytically using the basic 
properties of the classical digital computer (our main results are presented 
in [4, 5]). From our point of view, there are two general features of the 
computer when dealing with numbers: 1) numbers are represented as segments 
of a power series; 2) there is a digit carry procedure. These 
properties can be applied in the iterative or finite-difference scheme, 
where in the recurrence formulas the intermediate levels can be excluded. 
These properties allow us, for example, to represent solutions of differential 
equations (or systems of differential equations). At each level, we present the 
solution as the segment of a series in the powers of the step $\tau$ of the 
independent variable (we focus on the finite-difference scheme approximating 
the differential equation under consideration). Thus one can construct an 
analytical procedure, i.e. "computing without a computer", if using the 
computer analogy makes possible to exclude the intermediate expression (as a 
computer, in fact, excludes the intermediate values when operating with numbers). 
The construction of the proposed method does not only reduce the number of the 
arithmetical operations in calculations. It provides a solution in the 
explicit form (as a computer provides a solution in the numerical form, 
after executing many intermediate and "hidden" operations). The final analytic 
solution is treated as the limit when the small parameter $\tau$ tends to zero. 
In this short paper we describe the general features of the proposed method.

\section{The $\tau$-computer model}
In the present paper we consider general characteristics of our approach which 
can be applied when calculations are needed. A new theoretical 
model of the digital computer can be proposed if there is a small parameter in the problem under 
consideration. In this case one can use the expansion 
of the unknown function in a series in powers of this parameter. The universal 
approach of the finite-difference schemes (for solving differential equations) 
possesses such a parameter, namely a step on an independent variable. 

Consider the following Cauchy problem for the following ordinary differential equation: 
\begin{equation}
\frac {dy}{dt} = f\left(y\right),\ y\left(0\right)=y_0,
\label{eq:one}
\end{equation}
where $f\left(y\right)$ is an infinitely differentiable function with bounded derivatives. 
We will refer to the argument $t$ as "time". Consider a finite-difference explicit 
scheme of the following form:
\begin{equation}
y_{n+1} = y_n + \tau G\left(y_n\right),
\label{eq:two}
\end{equation}
where $G$ is also infinitely differentiable and is defined by the chosen numerical 
method. We assume that $\tau$ is small enough that scheme (2) is stable; 
hence, it converges to the solution of problem (1). 
One can see that after performing procedure (2) we obtain a segment of the power 
series (with powers of $\tau$), but the number of terms increases with $n$. We will search 
for a representation of the solution by this series with restricted number of terms. 
The digit carry procedure is introduced that guarantees that the absolute values 
of the coefficients of the segments of the series are less than $1/\tau$. If the coefficient 
becomes greater or equal to $1/\tau$, its part is carried to the more significant digit (i.e. left digit). 
If the absolute value of the coefficient is less than $1/\tau$, we call it "normalized"; 
if not, then we call it "not normalized" and denote it with the tilde $\tilde{a}_{i,n}$.
Consider the following representation of $y_n$:
\begin{equation}
y_n = \sum\limits_{i=0}^p \beta_i a_{i,n} \tau^i,
\end{equation}
where the factors $\beta_i \in \{-1,1\}$ are chosen such that all $a_i$ are 
not negative. This segment of the series provides a unique representation for the 
numbers with error $O\left(\tau^p\right)$. Suppose that we apply (2) to obtain $y_n$ 
from $y_{n-1}$. We cannot guarantee that for any $i$ the inequality 
$a_{i,n}<1/\tau$ holds true. Thus, 
\begin{eqnarray*}
\tilde{y}_n = \sum\limits_{i=0}^p \beta_i \tilde{a}_{i,n} \tau^i\ ,
\end{eqnarray*}
Next, we apply carry procedure starting from $\tilde{a}_{p,n}$:
\begin{eqnarray*}
&& a_{p,n} = \tilde{a}_{p,n} \bmod {\frac{1}{\tau}},\\  
&& \delta_{p,n} = \left(\tilde{a}_{p,n} - a_{p,n}\right)\tau = \tilde{a}_{p,n}\tau-\tilde{a}_{p,n}\tau \bmod {1} = \left[\tilde{a}_{p,n}\tau\right].
\end{eqnarray*}
For other coefficients, we obtain the following:
\begin{eqnarray*}
&& a_{m,n} = \left(\tilde{a}_{m,n} + \beta_m\beta_{m+1}\delta_{m+1,n}\right) \bmod {\frac {1}{\tau}},\\ 
&& \delta_{m,n} = \left[\left(\tilde{a}_{m,n} + \beta_m\beta_{m+1}\delta_{m+1,n}\right)\tau\right]\ .
\end{eqnarray*}
The condition $0\le a_i < 1/\tau$ implies that in the segments of the monotony, all of the 
coefficients $a_i$ increase monotonically; thus $a_{i,n} \le a_{i,n+1}$. Consider function $G$ in (2). 
Let us expand it into the Taylor series near $y^\ast=\beta_0 a_{0,n}$:
\begin{eqnarray*}
G\left(y_n\right)=G\left(y^\ast\right)+\sum\limits_{1}^{\infty} \left . {\frac {1}{k!}
\left(\sum\limits_{m=1}^{p} \beta_m a_{m,n}\tau^m\right)
\frac {d^k G\left(y\right)} {dy^k}} \right |_{y=y^\ast}.
\end{eqnarray*}
Using (1.2) and (1.3), we obtain the following:
\begin{eqnarray*}
y_{n+1} = \sum \limits _{m=0}^{p} \beta_m a_{m,n} \tau^m + \tau G\left(y^\ast\right) + 
\tau\sum\limits_{1}^{\infty} \left . {\frac {1}{k!}
\left(\sum\limits_{m=1}^{p} \beta_m a_{m,n}\tau^m\right)
\frac {d^k G\left(y\right)} {dy^k}} \right |_{y=y^\ast}.
\end{eqnarray*}
Thus a theoretical model of the proposed $\tau$-computer includes constructing an 
arbitrary positional numeral system, that is a 
mathematical auxiliary object that is useful in our approach; it deals with 
powers of the step of the argument. In contrast to the ordinary definition 
with the integer radix we consider the radix as the real number, namely, $1/\tau$. 
Consequently, the $\tau$-computer is an imaginary device that represents the numbers 
in the $\tau$-numeral system and applies carry procedure to prevent 
overflow of the digits. 

\section{The method of computer analogy}
The model of the theoretical computer will be completed if we will introduce 
important properties of the stochastic behavior of the less significant digits. It
allows us to reduce intermediate operations. Namely, the sums of stochastic 
numbers will be changed by appropriate means (expected values) which will be 
theoretically evaluated. In other words, "the inverse Monte Carlo procedure" 
will be realized. Then we will have the full computer model and the approach 
called the "method of computer analogy" (CA). 

The numbers that are produced by a certain formula are often referred to as 
pseudo-random numbers. These numbers are widely used in a variety of 
simulations and in Monte Carlo methods in general. An ideal random number 
generator (RNG) should produce numbers with the expected value $E\left[x\right]=1/2$ 
and the variance $Var\left(x\right) = 1/12$. 
Consider one of the most popular modern algorithms for generating random 
numbers by a linear congruent RNG, see e.g. [6]: 
\begin {eqnarray*}
x_m = \left(bx_{m-1} + c\right) \bmod{P},
\end {eqnarray*}
where $b$, $c$ and $P$ are integers and $b, c < P$. There are also more complicated 
formulae. Let us consider as a sample the following nonlinear system 
(and the initial problem): 
\begin {equation}
\left\{
\begin {aligned}
\frac {du}{dt} = v^2 - u^2, \ \ u\left(0\right)=1 \\
\frac {dv}{dt} = u^2 - 2v, \ \ v\left(0\right)=0.
\end {aligned}
\right .
\end {equation}
which has no an exact solution (in the ordinary sense). The finite-difference 
method (2) if we use the first-order scheme gives the following: 
\begin {eqnarray*}
\left\{
\begin {aligned}
u_{n+1} = u_n + \tau v_n^2 - \tau u_n^2, \ \ u_0 = 1, \\
v_{n+1} = v_n + \tau u_n^2 - 2 \tau v_n, \ \ v_0 = 0. \\
\end {aligned}
\right .
\end {eqnarray*}
According to CA we search for a solution in the following form: 
\begin {eqnarray*}
u_n = 1 - a_{1,n}\tau + a_{2,n}\tau^2 - a_{3,n}\tau^3, \ \ v_n = v_{1,n}\tau-b_{2,n}\tau^2+b_{3,n}\tau^3.
\end {eqnarray*}
It can 
be proved that for the quadratic nonlinearity it is sufficient to retain the 
third power in the segments of the series. It can also be proved that in the problem under consideration there is no 
carry to the zeros digits therefore the first (linear term) digits have the main 
roles (see [4,5]). The carry procedure gives the following: 
\begin {eqnarray*}
&& u_{n+1} = 1 - 
\left ( a_{1,n} + 1 - \delta_{2,n} \right ) \tau + 
\left ( \left ( a_{2,n} + 2 a_{1,n} - \delta_{3,n} \right ) \bmod {\frac{1}{\tau}} \right ) \tau^2 - \\ 
&& - \left ( \left ( a_{3,n} - b_{1,n}^2 + a_{1,n}^2 + 2a_{2,n} \right ) \bmod {\frac {1}{\tau}} \right )\tau^3\ , \\
&& v_{n+1} = 
\left ( b_{1,n} + 1 - \omega_{2,n} \right ) \tau - 
\left ( \left ( b_{2,n} + 2a_{1,n} + 2b_{1,n} - \omega_{3,n} \right ) \bmod {\frac {1}{\tau}} \right ) \tau^2 + \\
&& \left ( \left ( b_{3,n} + a_{1,n}^2 + 2a_{2,n} + 2b_{2,n} \right ) \bmod {\frac {1}{\tau}} \right ) \tau^3\ .
\end {eqnarray*}
where
\begin {eqnarray*}
&& \delta_{2,n} = \left [ \left ( a_{2,n} + 2 a_{1,n} - \delta_{3,n} \right ) \tau \right ],\ \ 
\omega_{2,n} = \left [ \left ( b_{2,n} + 2a_{1,n} + 2b_{1,n} - \omega_{3,n} \right ) \tau \right ]\ , \\
&& \delta_{3,n} = \left [ \left ( a_{3,n} - b_{1,n}^2 + a_{1,n}^2 + 2a_{2,n} \right ) \tau \right ],\ \  
\omega_{3,n} = \left [ \left ( b_{3,n} + a_{1,n}^2 + 2a_{2,n} + 2b_{2,n} \right ) \tau \right ]. 
\end {eqnarray*}
For simplicity denote $a = a_1$ and $b = b_1$. From this follows the expressions for the linear coefficients:
\begin {eqnarray*}
a_{n+1} = \sum \limits_{m=0}^{n}\left ( 1 - \delta_{2,m} \right ), \ \ b_{n+1} = \sum \limits_{m=0}^{n}\left ( 1 - \omega_{2,m} \right ). 
\end {eqnarray*}
Here we suppose that $\delta_2$ and $\omega_2$ are random integer values 
with the probabilities governed by the coefficients $a$ and $b$. 

This formula produces integer numbers in the segment $\left[0,\ 1/\tau-1\right]$. Fig. 1 shows the behaviour 
of coefficients that depend on the number of the level (the discretized time). 
Coefficients $a_{i,n}$ and $b_{i,n}$ depend on $a_{i,n-1}$ and $b_{i,n-1}$ respectively, and they can also indirectly 
(due to carries) depend on the coefficients of the previous level. Coefficients $a_{1,n}$ 
and $b_{1,n}$ do not exhibit stochastic behaviour, but they influence on all other values of 
the solution. Fig. 1 demonstrates the stochastic behaviour of coefficients $a_{i,n}$ and $b_{i,n}$. The stochastic 
behaviour of $a_{2,n}$ and $b_{2,n}$ slowly transform into regularity with oscillations. Nevertheless, 
we can still apply probabilistic methods when treating the carried values as 
weak dependent random values. 

\begin{figure} [ht]
\centering
\includegraphics{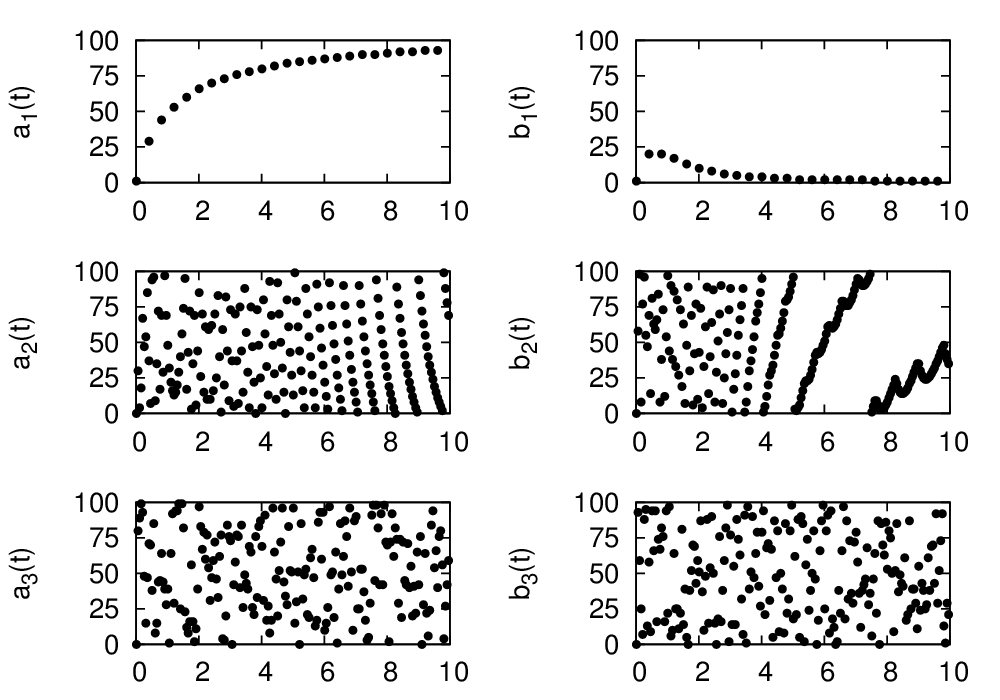}
\caption{Dependence of coefficients $a_1$, $a_2$, $a_3$, $b_1$, $b_2$, $b_3$ on time.}
\label{the-label-for-cross-referencing}
\end{figure}

The stochastic properties allow us to obtain the inverse 
dependence of the linear term on the discretized time. It appears that the 
coefficients of $\tau$ in the nonlinear powers have a stochastic behaviour; thus, we can 
apply the statistical methods. The carried value $\delta_{2,n}$ and $\omega_{2,n}$ are the only values 
that directly define the linear term. 
The next stage is expressing non-random terms using the carried values. 
Instead of carried values, we use their expected values. We refer to these values 
as quasi-random numbers. Let us use the law of large numbers to average the sequence 
of these numbers using their expected values. The $n^{-1}$ factor in the convergence 
formula for low discrepancy points may be contrasted with $n^{-1/2}$ convergence (see [7]). 
So the analog of the law of large numbers gives the following inequality: 
\begin {eqnarray*}
\left|\frac {1}{n} \sum\limits_{m=1}^{n}\delta_{2,m} - \frac{1}{n}\sum\limits_{m=1}^{n}E\left(\delta_{2,m}\right)\right|<\frac {c}{n},\ \ c=const.
\end {eqnarray*}
Here $E\left(\delta\right)$ denotes the expected value of a random variable $\delta$. It is admissible to 
assume that convergence is ordinary rather than in probability. The numeric simulation 
confirms that the error is close to $n^{-1}$. In Fig 2, one can see that the error does not 
exceed $n^{-1}$. plus some constant proportional to $\tau$. Thus the error can be estimated as $c/n+O\left(\tau\right)$. 
Since $\tau$ is of the order of $n^{-1}$ we conclude that resulting error is of the order of $n^{-1}$. 

\begin{figure} [ht]
\centering
\includegraphics{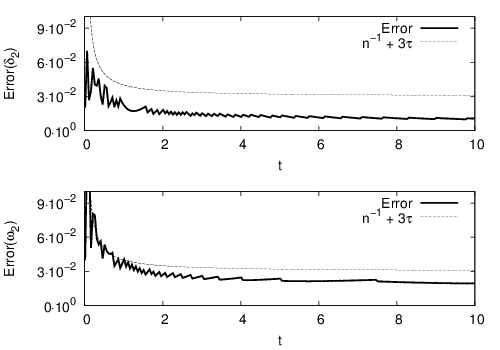}
\caption{Stochastic errors in CA.}
\label{the-label-for-cross-referencing}
\end{figure}

\section{Obtaining the explicit analytical form of the solution}
It is desirable to use the solutions without interpolation, therefore time intervals 
for $u$ and $v$ must be equal. For this purpose we use $a$ or $b$ as a basis variable, 
then the other variable is scaled proportionally to the ratio of time intervals, because 
for a small time interval $a$ and $b$ are changed linearly. We will use $a$ as a basis
variable. 

Let us denote the expected values of $1-\delta_2$ and $1-\omega_2$ as $E_\delta\left(a,b\right)$ and $E_\omega\left(a,b\right)$. 
The time intervals 
where $a$ and $b$ are not changed are $1/\left|E_\delta\left(a,b\right)\right|$ and $1/\left|E_\omega\left(a,b\right)\right|$ respectively.
The expected values are 
\begin {eqnarray*}
E_\delta\left(a, b\right) = \left(1-a\tau\right)^2 -b^2\tau^2, 
E_\omega\left(a, b\right) = \left(1-a\tau\right)^2 - 2b\tau.
\end {eqnarray*}
From this follows that the expected values of $\delta_2$ and $\omega_2$ are constant 
on layers defined by $a$ and $b$. Denoting the number of terms on $a$-th 
layer as $\Delta n_a = n_{a+1} - n_a$, we obtain:
\begin{eqnarray*}
a = \sum \limits_{m=1}^{n_a} \left(1-\delta_{2,n}\right) = 
\sum \limits_{m=1}^{n_a} E\left(1-\delta_{2,n}\right). 
\end{eqnarray*}
Since the increment $\Delta a$ equals 1, it follows that
\begin{eqnarray*}
E_\delta\left(m\right) \Delta n_m = 1 \Rightarrow n_{m+1} = n_m + 
\frac {1}{E_\delta\left(m\right)} = \sum \limits_{i=1}^m \frac {1}{E_\delta\left(i\right)}.
\end{eqnarray*}
As we have mentioned, we use $a$ as an independent variable, thus $n$ (the number of 
a level) and $b$ must be expressed through $a$:
\begin {eqnarray*}
&& u_a = 1 - a\tau, \ \ 
v_a = b_a\tau = \tau + \tau \sum\limits_{m=1}^{a-1} \prod\limits_{k=m+1}^{a-1}\left(1-\frac{2\tau}{\left(1-k\tau\right)^2}\right), \\
&&t_a = \tau n_a = \tau \sum \limits_{m=0}^{a-1}\frac{1} {\left(1-m\tau\right)^2}.
\end {eqnarray*}
The explicit solution allows us to construct approximations for asymptotic cases. 
Such asymptotic approximation valid for small (and moderate) $t$ is given in [5]. 
Namely 
\begin {eqnarray*}
&& u_a = 1 - a\tau, \ \ 
v_a = b_a\tau = \tau + \tau\left(1-2\tau\left(a+1\right)\right)\sum\limits_{m=2}^a\frac{1}{1-2m\tau} , \\
&& t_a = \tau n_a = \tau \sum \limits_{m=0}^{a-1}\frac{1} {1-2m\tau}.
\end {eqnarray*}
The important issue is finding the limit of this solution when the step size tends 
to zero, it will provide the final form of the solution. The convergence of the 
solutions at different values of the step sizes is shown in Fig. 3 (in the top the 
scaled segment of $v\left(t\right)$ is shown). One can see 
that the accuracy of the method is of the order of t. 
\begin{figure} [ht]
\centering
\includegraphics{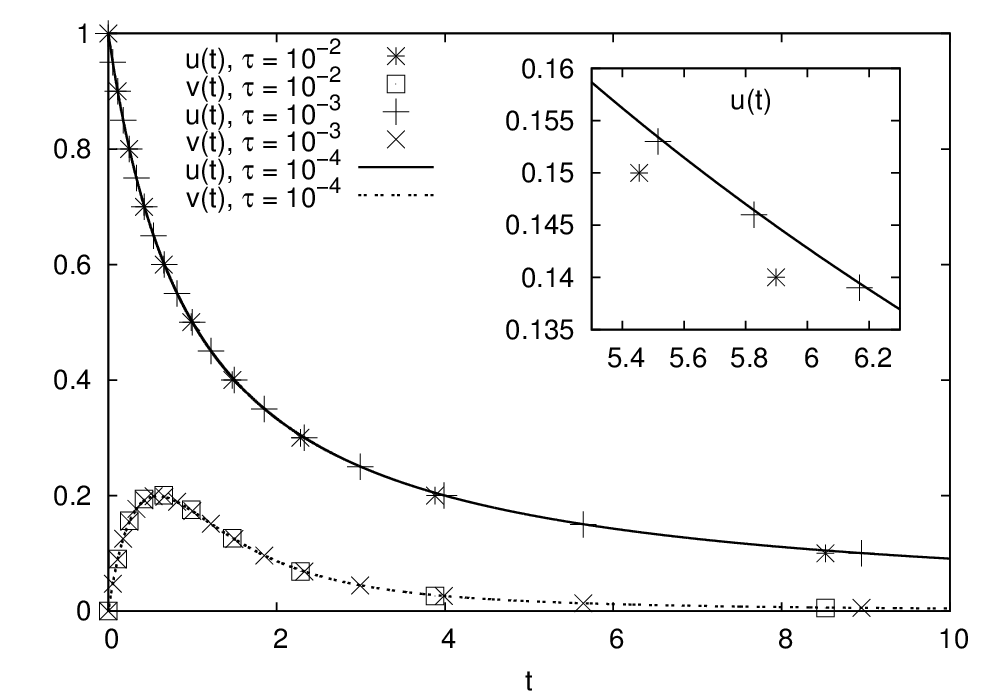}
\caption{Convergence to the solution.}
\label{the-label-for-cross-referencing}
\end{figure}

The passage to the limit 
implies that the step size $\tau$ tends to zero. Thus we obtain a limit form of the solution, 
which is a solution for the given problem in a certain range of the argument. 
Note, that the Sergeev computer [8], dealing with the infinitely small numbers, 
could be considered as a hypothetical tool for  
passing to the limit. 
\section{Concluding remarks}
In the present paper, we described the main features of a novel approach, that 
is based on the idea of a connection between analytical and computational mathematics. 
A theoretical model of the computer was presented. 
We analyzed the possibility of 
formalizing the computer operations, and proposed a representation 
of a solution in the form of a segment of the series in the powers of the step of 
the independent variable. This technique involves analytical work for obtaining the 
probabilities (which can be a difficult task). However, the main steps can be listed 
as follows: choose a convergent finite-difference scheme, find the number of terms to 
be retained, and obtain the expressions for the main coefficients in the representation of 
the solution. 

The proposed method of the computer analogy allows us to represent the 
solution of the problem as a convergent series that can be analyzed.

%
%

\end{document}